\documentclass[twoside,leqno,amsmath,thmmarks]{article}
\usepackage{pstricks,%
  pst-plot,%
  geometry,
  amsmath,%
  ntheorem,%
  quoting,%
}
\usepackage[author-year]{amsrefs}
\usepackage{lmodern}
\usepackage[T1]{fontenc}

\usepackage[ntheorem]{empheq}

\quotingsetup{leftmargin=2em, rightmargin=0in}

\usepackage{xspace}

\makeatletter
\DeclareRobustCommand\yspace{\futurelet\@let@token\@yspace}
\def\@yspace{%
  \ifx\@let@token\bgroup\else
  \ifx\@let@token\egroup\else
  \ifx\@let@token\/\else
  \ifx\@let@token\ \else
  \ifx\@let@token~\else
  \ifx\@let@token.\else
  \ifx\@let@token!\else
  \ifx\@let@token,\else
  \ifx\@let@token:\else
  \ifx\@let@token;\else
  \ifx\@let@token?\else
  \ifx\@let@token/\else
  \ifx\@let@token'\else
  \ifx\@let@token(\else
  \ifx\@let@token)\else
  \ifx\@let@token-\else
  \ifx\@let@token\@xobeysp\else
  \ifx\@let@token\space\else
  \ifx\@let@token\@sptoken\else
   \space
   \fi\fi\fi\fi\fi\fi\fi\fi\fi\fi\fi\fi\fi\fi\fi\fi\fi\fi\fi}
\makeatother

\makeatletter
\@ifundefined{usingAMS}
{
  \RequirePackage{amsmath}
  \RequirePackage{ntheorem}

  \newcommand*{\qedbox}{\lower.2ex\hbox{\vbox{\hrule
        \hbox{\vrule height 1.4ex
          \kern 1.4ex
          \vrule height 1.4ex}
        \hrule}}}

  \theorembodyfont{\itshape}
  \theoremstyle{margin}  

  \@ifundefined{PlainTheoremNumbering}
  {\@ifundefined{chapter}
    {\newtheorem{Theorem}{Theorem}[section]}
    {\newtheorem{Theorem}{Theorem}[chapter]}}
  {\newtheorem{Theorem}{Theorem}}

  \newtheorem{Lemma}[Theorem]{Lemma}
  \newtheorem{Corollary}[Theorem]{Corollary}

  \theoremsymbol{\qedbox}

  \theoremsymbol{}  

  \theoremstyle{plain} 
  \theoremheaderfont{\upshape\bfseries}
  \theorembodyfont{\upshape}
  \theoremsymbol{$\Diamond$}
  
  \theoremsymbol{}  

  \theoremstyle{nonumberplain}

  \theorembodyfont{\itshape}

  \theorembodyfont{\itshape}
  \newtheorem{Definition}{Definition}
  
  \newtheorem{Remark}{Remark}

  \theoremheaderfont{\upshape\bfseries}
  \theorembodyfont{\upshape}
  \theoremsymbol{$\Diamond$}
  
  \theoremsymbol{}  

  \theoremheaderfont{\scshape}
  \theorembodyfont{\upshape}
  \theoremseparator{.}
  \theoremsymbol{\qedbox}
  \newtheorem{Proof}{Proof}

  \newtheorem{Proof_of_the_theorem}{Proof of the theorem}

  \theoremstyle{nonumberbreak}
  
}
{\relax}
\makeatother

\newcommand*{\bfOregR}{\textbf{OR}\yspace}
\newcommand*{\OregR}{\ensuremath{\mathit{OR}}\yspace}

\def\set[[#1]]{\ensuremath{\left\{#1\right\}}\xspace}

\newcommand*{\dmu}{\,d\mu}
\newcommand*{\abs}[1]{\ensuremath{\left\lvert#1\right\rvert}} 
\newcommand*{\Nu}[1]{\ensuremath{\left\lVert#1\right\rVert}} 
\newcommand*{\mtbold}[1]{\ensuremath{\mathbf{#1}}}
\newcommand*{\R}{\mtbold{R}\yspace}  

\newcommand*{\figstart}{\hrule\bigskip}
\newcommand*{\figend}{\bigskip\hrule}

\newcommand*{\OR}{\textbf{OR}\yspace}

\newcommand*{\ds}{\,ds}
\newcommand*{\dt}{\,dt}

\newcommand*{\defeq}{\mathrel{\mathop{=}\limits^{\text{def}}}}

\title{A Best Possible General Form of the ``Master Theorem'' for Divide-and-Conquer Recurrences}

\author{Carl D. Offner\\
  University of Massachusetts Boston\\
  offner@cs.umb.edu}
\date{\today}

\begin{document}
\maketitle

\begin{abstract}
  \noindent
  We give here a general, best-possible, and smoothly-derived form of the
  ``Master Theorem'' for divide-and-conquer recurrences.
\end{abstract}

{\setlength{\parskip}{0pt}

  \tableofcontents
}

\section{Introduction}

The first 5 sections of this paper are expository, in preparation for what
follows, which is a new and best-possible result.

Section~\ref{sec:problem_and_history} states the problem and discusses the
various approaches that have been proposed, from the original Bentley-Haken-Saxe
Theorem  in \ocite{BentleyHakenSaxe1980} (later popularized as the ``Master
Theorem'' in \ocite{CormenLeisersonRivestStein2009}), to the greatly improved
version of \ocite{AkraBazzi1998}, which was further improved by
\ocite{Leighton1996}.

Section~\ref{sec:Akra_Bazzi_kernel} gives a smooth modern proof of the basic
result, based on some preparatory definitions and rewordings in
Sections~\ref{sec:functions_of_O_regular_variation} and
\ref{sec:back_to_the_Master_Theorem}.

Sections~\ref{sec:floors_and_ceilings} and \ref{sec:key_theorem} then give a
new, straightforward, and best possible way of dealing with floors and ceilings
in extending that basic theorem.

\section{The problem, and its history}
\label{sec:problem_and_history}

Algorithms falling into the ``divide-and-conquer'' category are typically
characterized by a simple recursion for their cost function $T(n)$ (where $n$ is
some measure of the size of the problem).  For instance, merge sort divides a
set of size $n$ into two sets of size $n/2$ (well, let's assume for the moment
that $n$ is even), then calls itself recursively on each of those sets and
finally merges the result.  It thus has a running time cost $T(n)$ that
satisfies
\begin{equation}
  \label{merge_sort_recurrence}
  T(n) = 2T\Bigl(\frac{n}{2}\Bigr) + cn
\end{equation}
where $c$ is some constant, and $cn$ is the cost of merging the two sorted lists
into one sorted list of size $n$.  (It's easy to see that the cost of this merge
is linear in $n$.)  It's not hard to see that any $T$ satisfying
(\ref{merge_sort_recurrence}) will grow as
\begin{equation*}
  T(n) = \Theta(n\log n)
\end{equation*}

For another example, suppose we consider the problem of multiplying two square
$n\times n$ matrices.  The Strassen algorithm is modelled by the recursion
\begin{equation*}
  T(n) = 7T\Bigl(\frac{n}{2}\Bigr) + cn^2
\end{equation*}
The key fact about this recursion is the appearance of 7 rather than 8 which
would appear in a more naive divide-and-conquer algorithm.  Because of the 7,
this recursion has solutions with growth
\begin{equation*}
  T(n) = \Theta(n^{\log_2 7})
\end{equation*}

A significant number of problems lead to recurrences of this kind, and while
each such recurrence could be solved\footnote{By ``solved'' we really mean that
  we are able to find tight bounds on the growth of the solution.  That's all we
  really care about, and exact closed-form solutions---i.e., in terms of
  elementary functions---are impossible in general anyway.}, it must have seemed
plausible that there would be a larger theory hovering in the background.

In any case, in 1980, Bentley, Haken, and Saxe published a short and charmingly
written paper that gave a method of deducing bounds for recursions of the form
\begin{equation*}
  T(n) = kT\Bigl(\frac{n}{c}\Bigr) + f(n)
\end{equation*}
Their method was based on the easily understood and intuitive notion of a
recursion tree which models the history of the execution of the algorithm.

The results of that paper, after some polishing, appeared in 1990 as the
``Master Theorem'' in the first edition of the widely used textbook
\ocite{CormenLeisersonRivestStein2009}.  Their theorem takes the following form:
\begin{Theorem}[``Master Theorem'']
  Let $a\geq 1$ and $b > 1$ be constants, let $g(n)$ be an ultimately positive
  function, and let $T(n)$ be defined on the non-negative integers by the
  recurrence
  \begin{equation*}
    T(n) = aT(n/b) + g(n)
  \end{equation*}
  where we interpret $n/b$ to mean either $\lfloor n/b\rfloor$ or $\lceil
  n/b\rceil$.  Then $T(n)$ has the following asymptotic bounds:
  \begin{enumerate}
  \item\label{MT:case_1} If $g(n) = O(n^{\log_b a - \epsilon})$ for some
    constant $\epsilon > 0$, then $T(n) = \Theta(n^{\log_b a})$.
  \item\label{MT:case_2} If $g(n) = \Theta(n^{\log_b a})$, then $T(n) =
    \Theta(n^{\log_b a}\log_2 n)$.
  \item\label{MT:case_3} If $g(n) = \Omega(n^{\log_b a + \epsilon})$ for some
    constant $\epsilon > 0$, and if
    \begin{equation}
      \label{weak_polynomial_growth}
      ag(n/b) \leq cg(n)
    \end{equation}
    for some constant $c < 1$ and all sufficiently large $n$, then $T(n) =
    \Theta\bigl(g(n)\bigr)$.
  \end{enumerate}
\end{Theorem}

It is actually more convenient\footnote{This seems to have been noticed first by
  \ocite{AkraBazzi1998}.} to work with functions defined on the reals, not
just on the integers.  So without belaboring the point, we will rewrite
our recursion as
\begin{equation*}
  T(x) = aT(x/b) + g(x)
\end{equation*}

\section{Functions of O-regular variation: basic definition}
\label{sec:functions_of_O_regular_variation}

We need to say something about the growth condition
(\ref{weak_polynomial_growth}).  \ocite{Leighton1996} found it useful to replace
this by a slightly different condition, which he called a \textbf{polynomial
  growth condition}.  And in fact, the class of functions satisfying that
condition was already a well-studied class known as the class of functions of
\textbf{O-regular variation}, the standard reference for which is Chapter 2
(Further Karamata Theory) of \ocite{Bingham-Goldie-Teugels1987}.  So we
henceforth refer to this class of functions as the class \bfOregR; this is the
standard abbreviation for it.  Appendix A of this paper contains a list of
equivalent properties of functions in this class; here is the one we use as its
definition:

\begin{Definition}
  The non-negative function $g$ on $[x_0,\infty$) (where $x_0\geq 0$) is a
  function of \emph{\textbf{O-regular variation}} iff there are constants
  $0 < A \leq 1 \leq B < \infty$ and $c > 1$ such that for all $x > x_0$ and for
  all $t \geq x_0$ such that
  \begin{equation*}
    \frac{x}{c} \leq t \leq x
  \end{equation*}
  we have
  \begin{equation}
    \label{Oreg_var}
    Ag(x) \leq g(t) \leq Bg(x)
  \end{equation}
\end{Definition}

This condition is in some respect a bit stronger than
(\ref{weak_polynomial_growth}) and in another respect somewhat weaker.  But it
doesn't really differ very much, and it is always satisfied in practice.  In
fact, other than ruling out functions that can't be bounded above and below by
powers of $x$, it does very little---any function that is bounded above and
below by some powers of $x$ but that fails to satisfy (\ref{Oreg_var})
will necessarily have oscillations that are large with respect to the size of
the function.

It is easy to see that a non-negative function $g\in\OregR$ must be either
positive everywhere or identically 0.  Further, if the function is non-zero,
then it must be bounded away from 0 and $\infty$ on every bounded interval.  And
it's trivially true that a function bounded away from 0 and $\infty$ on some
interval is of O-regular variation \textit{on that interval}.

\subsection{Subsets of the class \OregR}
\label{subsec:subsets_of_OR}

For convenience in some simple reasoning, we will make this definition:

\begin{Definition}
  For constants $0 < A \leq 1 \leq B < \infty$ and $c > 1$, we define $S(A, B,
  c)$ to be the class of functions $f:[0,\infty) \to [0,\infty)$ such that if
  \begin{equation*}
    \frac{x}{c} \leq y \leq x
  \end{equation*}
  then
  \begin{equation*}
    Af(x) \leq f(y) \leq B(x)
  \end{equation*}
\end{Definition}

The union of the sets $S(A, B, c)$ is thus all of \OregR.  And in fact, these
sets nest nicely, in two ways.  Here's the first way:

\begin{Lemma}
  \label{lemma:trivial_inclusion_for_S}
  If $1 < c < d < \infty$, then
  \begin{equation*}
    S(A, B, d) \subseteq S(A, B, c)
  \end{equation*}
\end{Lemma}
\begin{Proof}
  This is trivial but we'll write it out for clarity:  if $1 < c < d < \infty$, then
  \begin{equation*}
    \frac{x}{c} \leq y \leq x \implies \frac{x}{d} \leq y \leq x
  \end{equation*}
  and so
  \begin{empheq}{equation*}
    \begin{split}
      f \in S(A, B, d)
      &\iff \biggl(\frac{x}{d} \leq y \leq x \implies Af(x) \leq f(y) \leq Bf(x)\biggr)\\
      &\implies \biggl(\frac{x}{c} \leq y \leq x \implies Af(x) \leq f(y) \leq Bf(x)\biggr)\\
      &\iff f \in S(A, B, c)
    \end{split}
  \end{empheq}
\end{Proof}

We can think of this as saying that moving $c$ to the right makes $S(A, B, c)$
smaller.

We can also show that we can move $c$ to the right but \emph{not} decrease the set
$S(A, B, c)$ provided we are willing to adjust $A$ and $B$ at the same time.

Precisely, the following simple lemma shows that if $f\in\OregR$ with constants
$A$, $B$, and $c$, then we can stretch the interval
\begin{equation*}
  \frac{x}{c} \leq t \leq x
\end{equation*}
to the left, by increasing $c$, provided we are willing to also stretch the
interval $[A,B]$ in both directions.  And this is never a problem in practice.
\begin{Lemma}
  \label{lemma:stretching_interval}
  \begin{equation*}
    S(A, B, c) \subseteq S(A^2, B^2, c^2)
  \end{equation*}
\end{Lemma}
\begin{Proof}
  If $f\in S(A, B, c)$  and if
  \begin{equation*}
    \frac{x}{c^2} \leq y \leq x
  \end{equation*}
  then there are two possibilities:
  \begin{enumerate}
  \item If $y\in [x/c, x]$, then $Af(x) \leq f(y) \leq B(x)$, and so $A^2f(x)
    \leq f(y) \leq B^2f(x)$.

  \item Otherwise, $y\in [x/c^2, x/c)$, so $cy \in [x/c, x]$, and so we have
    \begin{equation*}
      Af(x) \leq f(cy) \leq Bf(x)
    \end{equation*}
    which yields the two inequalities
    \begin{equation*}
      A^2f(x) \leq Af(cy) \leq ABf(x) \text{      and      } ABf(x) \leq Bf(cy)
      \leq B^2f(x)
    \end{equation*}
    so
    \begin{equation}
      \label{eq:outer}
      A^2f(x) \leq Af(cy) \leq Bf(cy) \leq B^2f(x)
    \end{equation}
    Now since trivially, $y \in [y,cy]$, we also have
    \begin{equation}
      \label{eq:inner}
      Af(cy) \leq f(y) \leq Bf(cy)
    \end{equation}
    so putting together (\ref{eq:outer}) with (\ref{eq:inner}), we see that
    \begin{equation*}
      \begin{split} 
        A^2f(x) \leq f(y) \leq B^2f(x)
        \end{split}
    \end{equation*}
  \end{enumerate}
\end{Proof}

So this lemma shows that we can increase $c$ (thereby decreasing $x/c$) if we
are willing to decrease $A$ and increase $B$ at the same time.

\section{Back to the Master Theorem}
\label{sec:back_to_the_Master_Theorem}

Now getting back to the Master Theorem itself, it's actually not hard to see
that \emph{something} like this is going to be true.  For instance, if we set
$g$ to be identically 0, then it's easy to see that the resulting equation $T(x)
= aT(x/b)$ can be solved by trying the functions $T(x) = \alpha x^p$ (where
$\alpha$ is constant).  This will be a solution provided $ab^{-p} = 1$, which is
of course equivalent to $p = \log_b a$.  The actual proof given in the textbook
however is derived fairly directly from the reasoning in the original paper of
\ocite{BentleyHakenSaxe1980}.

There are some immediate questions that are raised by this result:
\begin{enumerate}
\item There are somewhat more general recursions that it seems might have
  similar behavior.  For instance, a number of algorithms lead to recursions of
  the form
  \begin{equation}
    \label{simple_recursion_n_terms}
    T(x) = \sum_{i=1}^n a_iT\Bigl(\frac{x}{b_i}\Bigr) + g(x)
  \end{equation}
  where $n = 2$; but of course in such a case it is worth looking to see if it
  is no harder to consider general values of $n \geq 1$.

  In order to handle such more general recursions, we are going to modify our
  notation just a little:

  \begin{itemize}
  \item First, we will assume (as will always be the case) that all the values
    \set[[b_i]] are contained in a fixed finite interval $[m,M]$ where
    $1 < m \leq M < \infty$.

  \item Next, let us define the ``delta measure'' $\delta_z$ to be
    the measure of mass 1 concentrated at $z$, and let us define the measure $\mu$
    by
    \begin{equation*}
      \mu = \sum_{i=1}^n a_i\delta_{b_i}
    \end{equation*}
  \end{itemize}
  Then $\mu$ is supported on the interval $[m,M]$ and instead of writing
  \begin{equation*}
    \sum_{i=1}^n a_if(b_i) \text{\qquad we can simply write \qquad}
    \int_m^M f(t)\dmu(t)
  \end{equation*}
  Similarly, the recursive equation
  \begin{equation*}
    T(x) = \sum_{i=1}^n a_iT\Bigl(\frac{x}{b_i}\Bigr) + g(x)
    \text{ \qquad becomes \qquad }
    T(x) \int_m^M T\biggl(\frac{x}{t}\biggr)\dmu(t) + g(x)
  \end{equation*}
  \and so on.  Now actually, everything we're doing from here on in works for
  \emph{any} finite positive measure $\mu$ supported on $[m,M]$, but we're only
  really concerned with the particular measure $\mu$ just defined above.  And
  we're using $\mu$ only as a notational convenience---everything we write could
  easily be rewritten in terms of sums as in (\ref{simple_recursion_n_terms}).

  Starting in Section~\ref{sec:Akra_Bazzi_kernel} below, we will consistently
  use this integral representation.

  Note that the identity $\sum_{i=1}^n a_ib_i^{-p} = 1$ is equivalent to
  \begin{equation*}
    \int_m^M y^{-p}\dmu(y) = 1
  \end{equation*}

\item The Master Theorem as stated applies to functions satisfying one of three
  possible conditions.  However, even disregarding functions $g$ that are wildly
  behaved and so wouldn't be expected to provide useful results, there are two
  ``gaps'' in the kinds of functions $g$ that are handled by this theorem:
  \begin{itemize}
  \item There is a gap between the functions $g$ handled in case \ref{MT:case_1}
    and \ref{MT:case_2}.  For instance, we could have a function
    \begin{equation*}
      g(x) = \frac{x^{\log_b a}}{\log x}
    \end{equation*}
  \item Similarly, there is a gap between the functions handled by cases
    \ref{MT:case_2} and \ref{MT:case_3}.
  \end{itemize}

\item Finally---and this is not so much a matter of the theorem itself, but of
  its proof---the theorem is actually not so difficult to prove if we disregard
  floors and ceilings.  The proof that introducing floors and ceilings is
  ``benign'' (in the sense that it does not alter the essential results of the
  theorem) is both computationally messy and uninformative---one doesn't feel
  after working through it that one has really learned anything.  It's just one
  of those painful things that you read once and then forget about.  This is
  dissatisfying, and it would be nice to find a good way to deal with this
  problem that was natural and intuitive.
\end{enumerate}

These three issues have all been considered, and to some extent resolved.
Curiously enough, this all happened in 1996:
\begin{enumerate}
\item If we just consider solutions of the form $\alpha x^p$ to the simplified
  recursion
  \begin{equation*}
    T(x) = \sum_{i=1}^n a_iT\Bigl(\frac{x}{b_i}\Bigr)
  \end{equation*}
  then as before we see that we have solutions precisely when
  \begin{equation}
    \label{sum_formula_for_p}
    \sum_{i=1}^n a_i b_i^{-p} = 1
  \end{equation}
  Thus, it seems reasonable to believe that functions of the form $x^p$ play a
  role in solutions to (\ref{simple_recursion_n_terms}) similar to the role they
  play in the Master Theorem.

  While I don't really know, it seems likely to me that a number of people must
  have noticed this.  In any case, a theorem of this form appears in
  \ocite{Kao1997}\footnote{The work was actually presented in 1996; see the
    bibliographic reference.}; however, he restricts himself to functions $g$ of
  the form
  \begin{equation*}
    g(x) = c x^\alpha \log^\beta x
  \end{equation*}

\item The big advance in this area was made by \ocite{AkraBazzi1998}, who
  also\footnote{Their paper was actually submitted in 1996} arrived at
  (\ref{sum_formula_for_p}) and in addition---and this was the key result---gave a
  more precise formula for the growth of $T$ which avoids having three cases
  with ``gaps'' in between.

  Akra and Bazzi proved that provided
  \begin{itemize}
  \item all $a_i > 0$ and $\sum_{i=1}^n a_i \geq 1$; and
  \item all $b_i > 1$; and
  \item $g$ is a non-negative function in \bfOregR\footnote{In their paper, Akra
    and Bazzi assumed a slightly stronger condition, but this is all that is
    necessary.}
  \end{itemize}
  then any solution of the recurrence
  $T$ of (\ref{simple_recursion_n_terms}) will satisfy
  \begin{equation}
    T(x)
    = \Theta\Biggl(x^p\biggl(1 + \int_{x_0}^x t^{-p}g(t)\frac{dt}{t}\biggr)\Biggr)
  \end{equation}
  where $x_0$ is any fixed number $> 1$.

  It is easy to see that the three cases of the Master Theorem follow quickly
  from this result\footnote{Well actually, this method forces us to make our
    peace with a very slight increase in strength of the hypotheses of the
    theorem---we now must assume that $g\in\OregR$ in \emph{all} cases, and
    further, that condition is slightly different from the original condition of
    the theorem as stated in Case 3.  However, the condition $g\in\OregR$ is in
    practice always satisfied anyway.}

  This is really an elegant result\footnote{We should also note that in the same
    year, \ocite{Wang_Fu1996} produced a very similar result, but only for
    recursions of the form considered in the Master Theorem, although they did
    also consider the more general case in which the numbers $a$ and $b$ were
    well-behaved functions of $n$.}---and the kernel of Akra and Bazzi's
  original proof can be made to be very smooth, as we show below in
  Section~\ref{sec:Akra_Bazzi_kernel}.  However, the problem of introducing
  floors and ceilings was still something that they could only deal with by a
  lengthy and painful addendum.

\item Very shortly after Akra and Bazzi announced their result (and before the
  result actually appeared), \ocite{Leighton1996} came up with another way to
  deal with the ``floors and ceilings'' problem.

  Leighton assumed that $g$ is in the class \OregR.  He then replaced
  (\ref{simple_recursion_n_terms}) with the more general recursion
  \begin{equation}
    \label{general_recursion_n_terms}
    T(x) = \sum_{i=1}^n a_i T\Bigl(\frac{x}{b_i} + h_i(x)\Bigr) + g(x)
  \end{equation}
  Here the functions $h_i(x)$ can be defined so that
  \begin{equation*}
    \frac{x}{b_i} + h_i(x) = \begin{cases}
      \displaystyle \biggl\lceil\frac{x}{b_i}\biggr\rceil \\
      \phantom{\biggl\lceil}\text{or}\\
      \displaystyle \biggl\lfloor\frac{x}{b_i}\biggr\rfloor
    \end{cases}
  \end{equation*}
  Of course if we do this then $\abs{h_i(x)} < 1$ for all $i$.  However,
  Leighton noticed that actually weaker bounds on the functions $h_i(x)$ could
  be accommodated; in fact they could be unbounded provided they did not grow
  too fast.  He showed that the Akra-Bazzi results remained true for
  (\ref{general_recursion_n_terms}) provided that for some $\alpha > 1$ and all
  $i$, and sufficiently large values of $x$,
  \begin{equation*}
    \abs{h_i(x)} \leq \frac{x}{\log^\alpha x}
  \end{equation*}
  and he also noted that one could not have $\alpha = 1$ and get to the same
  conclusion.
  
  Again, this result is itself very nice.  It in turn leads to two questions,
  however:
  \begin{itemize}
  \item The constraint on $h_i(x)$ looks a bit unnatural.  It seems likely that
    there is a more natural condition waiting to be found.
  \item It remains true that Leighton's proof is terminally clever and
    forbiddingly computational: It involves a number of heroic estimates as well
    as an inductive hypothesis that is seemingly out-of-the-blue.
  \end{itemize}

  It would certainly be nice to come up with a simpler and more motivated line
  of reasoning.

  It turns out that both of these questions can be resolved quite nicely, and we
  will do so below.
\end{enumerate}

\section{The core of  the Akra-Bazzi argument}
\label{sec:Akra_Bazzi_kernel}

\begin{Theorem}[Akra-Bazzi] \label{thm:Akra-Bazzi}
  We assume the following conditions are satisfied:
  \begin{enumerate}
  \item There is a finite (positive) measure $\mu$ whose support is contained in
    $[m,M]$, where $1 < m \leq M < \infty$, and such that $\Nu{\mu} \geq 1$.

  \item There is a positive measurable function $g:\R^+\to \R^+$ with
    $g\in\OregR$,

  \item There is a positive measurable function $T:\R^+\to\R^+$ and a number $x_0
    > M$ such that
    \begin{enumerate}
    \item $T(x)$ is bounded away from 0 and $\infty$ on $[1,mx_0]$, and
    \item For all $x \geq x_0$,
      \begin{equation}
        \label{recursion_no_h}
        T(x) = \int_m^M T\left(\frac{x}{t}\right)\dmu(t) + g(x)
      \end{equation}
    \end{enumerate}
  \end{enumerate}
  Under these conditions,
  \begin{equation}
    \label{ab_no_h}
    T(x)
    = \Theta\Biggl(x^p\biggl(1 + \int_{x_0}^x t^{-p}g(t)\frac{dt}{t}\biggr)\Biggr)
  \end{equation}
  where $p$ is the unique non-negative solution to
  \begin{equation}
    \label{formula_for_p}
    \int_m^M t^{-p}\dmu(t) = 1
  \end{equation}
\end{Theorem}

Before proving this, there are some simple points we need to make:
\begin{itemize}
\item The constraint that $T$ is bounded away from 0 and $\infty$ on $[1,mx_0]$
  is not a very strong one.  We really only care about the behavior of $T(x)$
  for large $x$, and as we will see, we can pretty much modify it as we want to
  for small values of $x$.

  Similarly, it is clear from (\ref{recursion_no_h}) and (\ref{ab_no_h}) that we
  only care about the values of $g(x)$ for $x \geq x_0$.

\item We should show that $p$ exists as stated.  This is easy:  if we set
  \begin{equation*}
    \phi(p) = \int_m^M t^{-p}\dmu(t)
  \end{equation*}
  then $\phi$ is differentiable and its derivative is strictly less than 0 for
  $p > 0$.  So $\phi$ is strictly decreasing for all $p\geq 0$.  Further, by
  our assumptions, $\phi(0) \geq 1$ and clearly $\phi(p)\to 0$ as $p\to\infty$.
  Thus there is a unique $p\geq 0$ for which $\phi(p) = 1$.

  Note that if $g\equiv 0$, then the solutions to (\ref{recursion_no_h}) are
  just all the functions of the form $T(x) = cx^p$, and this is consistent with
  (\ref{ab_no_h}).

\item Figure~\ref{fig:main_variables_used_in_proof} depicts the main variables
  we are dealing with and shows how they are related.  It's useful in following
  the reasoning below.

  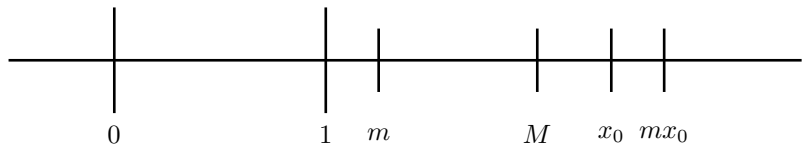
\begin{figure}[htbp]
    \figstart
    \begin{center}
      \psset{unit=0.4em,linewidth=1pt}
      \begin{pspicture}(0,-10)(75,10)
        \rput(0,0){\psline(0,0)(75,0)}
        \rput(10,0){\psline(0,-5)(0,5)}
        \rput(10,-7){0}
        \rput(30,0){\psline(0,-5)(0,5)}
        \rput(30,-7){1}
        \rput(35,0){\psline(0,-3)(0,3)}
        \rput(35,-7){$m$}
        \rput(50,0){\psline(0,-3)(0,3)}
        \rput(50,-7){$M$}
        \rput(57,0){\psline(0,-3)(0,3)}
        \rput(57,-7){$x_0$}
        \rput(62,0){\psline(0,-3)(0,3)}
        \rput(62,-7){$mx_0$}
      \end{pspicture}
    \end{center}
    \caption{The main variables used in the proof.}
    \label{fig:main_variables_used_in_proof}
    \figend
  \end{figure}

\end{itemize}

We need a preparatory lemma:
\begin{Lemma}
  \label{T_in_OR_lemma}
  Under conditions (1)--(3) of Theorem~\ref{thm:Akra-Bazzi}, $T\in\OregR$ on
  $[x_0, \infty)$.
\end{Lemma}
\begin{Proof}
  Since $T(x)$ is bounded away from 0 and $\infty$ on $[1,x_0m]$, where
  $x_0 \geq M$, we know at least that there are constants
  $0 < \Gamma_0 \leq 1 \leq \Gamma_1 < \infty$ such that for every $x$ and $y$
  in that interval,
  \begin{equation}
    \label{pre_OR_for_T}
    \Gamma_0T(x) \leq T(y) \leq \Gamma_1 T(x)
  \end{equation}
  We want to use this, together with the recursion(\ref{recursion_no_h}) to show
  that $T$ is in \OR on $[x_0,\infty)$.

  We define a sequence of intervals $\set[[I_j: 0 \leq j < \infty]]$ that
  partition $[x_0,\infty)$ as follows:
  \begin{equation*}
    I_j = [x_0m^j, x_0m^{j+1})
  \end{equation*}

  We will show that there is a constant $c_0 > 1$ and constants
  $0 < C_0 \leq 1 \leq C_1 < \infty$ such that for each $j \geq 0$ and for all
  $x\in I_j$,
  \begin{equation}
    \label{OR_for_T}
    C_0 T(x) \leq T(y) \leq C_1T(x) \quad\text{for}\quad  x/c_0 \leq y \leq x
  \end{equation}
  Let us first look at the interval $I_0$.  Let us set $c_0 = m$.  Since
  \begin{itemize}
  \item  $I_0 = [x_0, x_0m) \subset [1, x_0m]$
  \item and since (\ref{pre_OR_for_T}) holds for all $x$ and $y$ with
    $x\in [1,x_0m)$ and with $x/c_0 \leq y \leq x$,
  \end{itemize}
  we see that (\ref{OR_for_T}) holds for $I_0$ with $C_0 = \Gamma_0$ and
  $C_1 = \Gamma_1$.  ($c_0 = m$ is now fixed.  $C_0$ and $C_1$ may be adjusted
  once below, but otherwise will be fixed from here on.)

  Now since $g\in\OregR$ on $\R^+$, we know that $g$ satisfies a condition like
  this over $(0,\infty)$, although likely with different constants.  We will
  need to relate these to the constants used above for $T$.

  Let us therefore assume that $g\in S(A, B, c)$.  We know that $c > 1$.

  There are two cases we need to consider when relating $g$ to $T$:

  \begin{description}
  \item[Case 1: $c \geq c_0$.]\phantom{asdf}\linebreak
    In this case, by Lemma~\ref{lemma:trivial_inclusion_for_S}, $g\in S(A, B,
    c_0)$.
  \item[Case 2: $c < c_0$.]\phantom{asdf}\linebreak
    Take the smallest integer $n > 0$ such that $c^{2^n} \geq c_0$.  Then (by
    Lemma~\ref{lemma:stretching_interval}) we have
    $g \in S(A^{2^n}, B^{2^n}, c^{2^n})$, and then again by
    Lemma~\ref{lemma:trivial_inclusion_for_S}, $g \in S(A^{2^n}, B^{2^n}, c_0)$.
  \end{description}

  If necessary, we can then adjust $C_0$ downwards (still remaining > 0) and $C_1$
  upwards so that  $g \in S(C_0, C_1, c_0)$ as well.  And (\ref{OR_for_T}) will
  of course then continue to hold as well.

  At this point, $c_0$, $C_0$, and $C_1$ are fixed.  And (\ref{OR_for_T}) holds
  for the interval $[x_0m^j, x_0m^{(j+1)}]$ when $j = 0$.

  Now we proceed by induction.  We assume the result is true up through $j-1$;
  in particular, we assume that (\ref{OR_for_T}) holds for $x$ in
  $I_{j-1} = [x_0m^{j-1}, x_0m^j]$.  Then for
  \begin{itemize}
  \item $x_0m^j \leq x \leq x_0m^{j+1}$ and
  \item $\displaystyle\frac{x}{c_0} \leq y \leq x$ and
  \item $t\in[m,M]$
  \end{itemize}
  we have
  \begin{equation*}
    \frac{x}{t} \leq \frac{x}{m} \leq x_0m^j
    \end{equation*}
  and so, bearing in mind that $c_0 = m > 1$, we have
  \begin{alignat*}{2}
    \frac{1}{c_0}\frac{x}{t} &< \frac{y}{t} \leq \frac{x}{t} &\quad& \text{which
      we use immediately below  when estimating $T$} \\
    \frac{x}{c_0} &< y \leq x && \text{which we use immediately below when
      estimating $g$}
  \end{alignat*}
  and so we can see, using the recursion (\ref{recursion_no_h}), that
    
  \begin{align*}
    C_0T(x) &= \int_m^MC_0T\left(\frac{x}{t}\right)\dmu(t) + C_0g(x)\\
            &\leq \int_m^MT\left(\frac{y}{t}\right)\dmu(t) + g(y)
              \text{\qquad\qquad(and this is just $T(y)$)}\\
            &\leq\int_m^MC_1T\left(\frac{x}{t}\right)\dmu(t)  + C_1g(x)\\
            &=C_1T(x)
  \end{align*}
  which shows that (\ref{OR_for_T}) continues to hold for $x$ in the interval
  $I_j = [x_0m^j, x_0m^{(j+1)}]$; and this concludes the proof of the inductive
  step.  So $T\in\OregR$, and that in turn concludes the proof of the Lemma.
\end{Proof}

\begin{Corollary}
  Lemma~\ref{T_in_OR_lemma} holds with $c_0 = m$, and with some values of $C_0$
  and $C_1$.  And we may arrange things so that the same constants work for $g$.
\end{Corollary}
\begin{Proof}
  This follows from Lemma~\ref{lemma:stretching_interval}.
\end{Proof}

\begin{Proof_of_the_theorem}
  What Akra and Bazzi did---and this is the main contribution of their
  paper---is to derive (\ref{ab_no_h}), which shows that the growth of $T$ can
  be directly related to the growth of the integral
  \begin{equation*}
    \int_{x_0}^x t^{-p}g(t)\frac{\dt}{t}
  \end{equation*}

  Here is what they did:  From (\ref{recursion_no_h}) we have
  \begin{equation}
    \label{integrated_recursion}
    \int_{x_0}^x t^{-p}T(t)\frac{dt}{t}
    = \int_{x_0}^x t^{-p} \int_m^MT\Bigl(\frac{t}{s}\Bigr)\dmu(s)\frac{dt}{t}
    + \int_{x_0}^x t^{-p} g(t)\frac{dt}{t}
  \end{equation}

  The first term on the right of (\ref{integrated_recursion}) is

  \begin{equation*}
    \begin{aligned}
      \int_{x_0}^x t^{-p}\int_m^MT\Bigl(\frac{t}{s}\Bigr)\dmu(s)\frac{dt}{t}
      &= \int_m^M\int_{t=x_0}^x t^{-p}T\Bigl(\frac{t}{s}\Bigr)\frac{dt}{t}\dmu(s) \\
      &=\int_m^M s^{-p}\int_{t = \frac{x_0}{s}}^{\frac{x}{s}}
      t^{-p}T(t)\frac{dt}{t}\dmu(s) \\
    \end{aligned}
  \end{equation*}
  and the right-hand side of this then can be decomposed to into three terms, so
  we get the first term on the right of (\ref{integrated_recursion}) to be
      
  \begin{equation}
    \label{first_term_of_integrated_recursion}
    \begin{aligned}
      \int_{x_0}^x t^{-p}\int_m^MT\Bigl(\frac{t}{s}\Bigr)\dmu(s)\frac{dt}{t}
      &= \quad \int_m^M s^{-p} \int_{t=x_0}^x t^{-p}T(t)\frac{dt}{t}\dmu(s) \\
      &\quad - \int_m^M s^{-p} \int_{t=\frac{x}{s}}^x
      t^{-p}T(t)\frac{dt}{t}\dmu(s) \\
      &\quad
      + \int_m^M s^{-p} \int_{t=\frac{x_0}{s}}^{x_0} t^{-p}T(t)\frac{dt}{t}\dmu(s)
    \end{aligned}
  \end{equation}
      
  Now since $\int_m^M s^{-p}\dmu(s) = 1$, the first term on the right of
  (\ref{first_term_of_integrated_recursion}) becomes
  \begin{equation*}
    \int_m^M s^{-p} \biggl(\int_{t=x_0}^x t^{-p}T(t)\frac{dt}{t}\biggr)\dmu(s)
    = \int_{x_0}^x t^{-p} T(t)\frac{dt}{t}
  \end{equation*}
  which is just the left-hand side of (\ref{integrated_recursion}).
  Therefore, (\ref{integrated_recursion}) becomes
  \begin{equation}
    \label{zero_form}
    0 =     - \int_m^M s^{-p} \int_{t=\frac{x}{s}}^x
      t^{-p}T(t)\frac{dt}{t}\dmu(s) \\
      + \int_m^M s^{-p} \int_{t=\frac{x_0}{s}}^{x_0} t^{-p}T(t)\frac{dt}{t}\dmu(s)
      + \int_{x_0}^x t^{-p} g(t)\frac{dt}{t}
    \end{equation}

    \begin{itemize}
    \item Because in addition $T\in\OregR$, the first term on the right of (\ref{zero_form})
      can be bounded as follows:
    \begin{equation*}
      C_0x^{-p}T(x)\frac{m^p - 1}{p}
      \leq \int_m^Ms^{-p} \int_{t=\frac{x}{s}}^x t^{-p}T(t)\frac{dt}{t}\dmu(s)
      \leq C_1x^{-p}T(x)\frac{M^p - 1}{p}
    \end{equation*}

  \item And the second term (while it could be bounded similarly) is just a
    positive constant, which we will denote by $C$.
  \end{itemize}
  So (\ref{zero_form}) thus produces the bounds
  \begin{equation*}
    C_0x^{-p}T(x)\frac{m^p - 1}{p}
    \leq \int_{x_0}^x t^{-p}g(t)\frac{dt}{t} + C
    \leq C_1x^{-p}T(x)\frac{M^p - 1}{p}
  \end{equation*}
  which leads immediately to
  \begin{equation*}
    T(x)
    = \Theta\Biggl(x^p\biggl(1 + \int_{x_0}^x t^{-p}g(t)\frac{dt}{t}\biggr)\Biggr)
  \end{equation*}
  and that concludes the proof of Theorem~\ref{thm:Akra-Bazzi}.
\end{Proof_of_the_theorem}

\section{Dealing with floors and ceilings}
\label{sec:floors_and_ceilings}

\subsection{The functions $\lambda$ and $\mu$}
\label{subsec:lambda_and_mu}

We will follow Leighton in considering the modified recursion
\begin{equation}
  \label{general_recursive_formula}
  T(x) = \int_m^MT\Bigl(\frac{x}{t} + h(x,t)\Bigr)\dmu(t) + g(x)
\end{equation}

We know (under the conditions of the Akra-Bazzi theorem above) that when
$h(x,t)$ is identically 0, that
\begin{equation}
  \label{akra_bazzi_result}
  T(x)
  = \Theta\Biggl(x^p\biggl(1 + \int_{x_0}^x t^{-p}g(t)\frac{dt}{t}\biggr)\Biggr)
\end{equation}
We want to find conditions on the functions $h(\cdot,t)$ such that this remains true.

\begin{quoting}
Now there is one point just crying out to be made here:  it really looks like
the recursion (\ref{general_recursive_formula}) might be thought of as some sort
of perturbation of the same formula with $h(x,t)\equiv 0$.  It would be really
nice to be able to say---even without knowing the Akra-Bazzi theorem---that
given some natural-looking constraints on the functions $h(x,t)$, the solutions to
(\ref{general_recursive_formula}) of necessity behave the same as those in the
Akra-Bazzi theorem.

However, I have not been able to come up with such a theorem, and as far as I
know, no one else has either.  Therefore, we will in the following just prove
our result from scratch, by an induction.  This is somewhat similar to what
Leighton did, although we perform an induction using the multiplicative
structure of the reals rather than the additive structure.  This enables our
reasoning to make more natural use of the multiplicative convolution in
(\ref{recursion_no_h}); it is as a consequence more elementary, and the result
we arrive at is improved from previous results and in fact is best possible for
a very slightly restricted version of the problem, as described just below in
Section~\ref{subsec:restrict}.

And to start with, we will not actually use the Akra-Bazzi theorem---our result
will be proved ``from scratch''.
\end{quoting}

First of all we can see that the functions $h(x,t)/x$ must be uniformly bounded
above---in fact, by $1 - \frac{1}{m}$.  If this were not so, then
(\ref{general_recursive_formula}) would not serve to define $T(x)$ in terms of
values of $T(y)$ where $y < x$, and that is an essential point of what we are
doing.  Similarly, the functions $h(x,t)/x$ must be uniformly bounded below,
since $T(x)$ is defined only for $x\geq 0$.

Thus, we know that there are bounded functions $\lambda(x) \geq 0$ and
$\mu(x)\geq 0$ defined on $[1,\infty)$ such that for each $t \in [m,M]$,
\begin{equation*}
  -\lambda(x) \leq \frac{h(x,t)}{x} \leq \mu(x)
\end{equation*}

\subsection{We restrict the problem, very slightly}
\label{subsec:restrict}

At this point, we make a simplification: we will henceforth restrict our
attention to functions that ``don't wiggle too much''.  Precisely, we will
insist that $\lambda$ and $\mu$ are non-increasing on $[1,\infty)$.  Therefore
each of them decreases to a non-negative limit as $x\to\infty$.  One might in
principle try to relax this restriction, but it's hard to see how this added
generality would yield anything useful for the problem we are addressing here.

Then, continuing, let us define the class of functions
\begin{equation*}
  H(\lambda, \mu)
  = \set[[h: \forall x\geq x_0 \text{ and } t\in[m,M], -x\lambda(x) \leq h(x,t) \leq x\mu(x)]]
\end{equation*}

Our goal is to find conditions on $\lambda$ and $\mu$ which will ensure that if
$h\in H(\lambda,\mu)$ then all solutions of the corresponding recursion satisfy
the Akra-Bazzi bounds.

A necessary condition for this to hold is that both
$\lim_{t\to\infty}\lambda(t)$ and $\lim_{t\to\infty}\mu(t)$ must be 0.  For if
$\lim_{t\to\infty}\mu(t) = \epsilon > 0$, say, then if we define $h(x,t) =
\epsilon x$ for all $t$ and take $g\equiv 0$, we can write
\begin{equation*}
  \frac{x}{t} + \epsilon x = \frac{x}{\alpha(t)}
\end{equation*}
where
\begin{equation*}
  \alpha(t) = \frac{t}{1 + \epsilon t} < t
\end{equation*}
and so (\ref{general_recursive_formula}) with $g\equiv 0$ becomes the recursive formula
\begin{equation*}
  T(x) = \int_m^MT\Bigl(\frac{x}{\alpha(t)}\Bigr)\dmu(t)
\end{equation*}
This has solutions of the form $cx^q$, where $q$ is the unique solution of
\begin{equation*}
  \int_m^M \alpha(t)^{-q}\dmu(t) = 1
\end{equation*}
Now since $\alpha(t) < t$, we know that $q < p$ and this is incompatible with
(\ref{akra_bazzi_result}).  An entirely similar argument shows that
$\lim_{t\to\infty}\lambda(t) = 0$.

\subsection{A property of the function $g$}

We have already assumed that $g$ is a non-negative measurable function in the
class \OregR.  In everything that follows, $m$ and $M$ will be what they have
always been---their values in equation (\ref{formula_for_p}), and are fixed.  We
know that $1 < m \leq M < \infty$.

We know, therefore, that there are numbers $c > 1$  and $x_0 > 0$ as well as
constants $0 < A \leq 1 \leq B < \infty$ such that for all $x > x_0$ and for all
$t$ such that
\begin{equation*}
  \frac{x}{c} \leq t \leq x
\end{equation*}
we have
\begin{equation*}
  Ag(x) \leq g(t) \leq Bg(x)
\end{equation*}

Now in preparation for what follows, we are going to extend this notation a bit:
We know by the results of Section~\ref{subsec:subsets_of_OR} that we may take
$c$ as large as we want (adjusting $A$ and $B$ accordingly).

We can also take $x_0$ to be as large as necessary, as well.  And in any case,
we know already from the assumptions in Theorem~\ref{thm:Akra-Bazzi} that
$x_0 > M > 1$.

So we will introduce a new variable $m_0$ such that
\begin{equation*}
  1 < m_0 < m
\end{equation*}
And then we will make $c$ larger than $M$, and define $M_0 = c$, so that
\begin{equation*}
  M < M_0 < \infty
\end{equation*}
And finally, we will increase $x_0$ if necessary so that $M_0 < x_0$.  $x_0$ may
be increased once again in Section~\ref{subsec:props_of_h}, but is otherwise fixed.

Now with
\begin{equation}
  \label{eq:m0_and_M0}
  1 < m_0 < m \leq M < M_0  = c < \infty
\end{equation}
we have for all $t$ such that
\begin{equation*}
  \frac{x}{M_0} \leq t \leq x
\end{equation*}
and all $x > x_0$,
\begin{equation}
  \label{g_is_nice_past_x_0}
  Ag(x) \leq g(t) \leq Bg(x)
\end{equation}
As an immediate consequence we see that
\begin{Lemma}
  for each $m_0$ and $M_0$ as above, there must be constants
  $0 < C \leq D < \infty$ such that for all $x > M_0$,
  \begin{equation}
    \label{bounds_on_g_integral}
    C g(x)
    \leq x^p\int_{\frac{x}{m_0}}^xt^{-p}g(t)\frac{dt}{t}
    \leq x^p\int_{\frac{x}{M_0}}^xt^{-p}g(t)\frac{dt}{t}
    \leq D g(x)
  \end{equation}
\end{Lemma}
\begin{Proof}
  The middle inequality is automatically true because the integrand is
  non-negative.  The other two inequalities follow from
  (\ref{g_is_nice_past_x_0}).
\end{Proof}

\subsection{Properties of the functions $h(x,t)$}
\label{subsec:props_of_h}

Now let us see what we can say about the functions $h(x,t)$.

By what we have proved above we know that given $m_0$ and $M_0$ with
$1 < m_0 < m \leq M < M_0 < \infty$, we can increase $x_0$ if necessary (making
use of the fact that $\lambda(x)$ and $\mu(x)$ tend to 0 as $x\to\infty$), to
ensure that for all $x > x_0$ and for all $t\in[m,M]$,

\begin{equation}
  \label{eq1}
  \frac{x}{M_0}
  \leq x\Bigl(\frac{1}{M} - \lambda(x)\Bigr)
  \leq \frac{x}{t} + h(x,t)
  \leq x\Bigl(\frac{1}{m} + \mu(x)\Bigr)
  \leq \frac{x}{m_0}
\end{equation}
Also, and in any case (for all $x > x_0$ and for all $t\in[m,M]$), we must have
\begin{equation}
  \label{eq2}
  1 - M\lambda(x)
  \leq 1 + t\frac{h(x,t)}{x}
  \leq 1 + M\mu(x)
\end{equation}

We will assume that the numbers $m_0$ and $M_0$ are fixed---their exact values
don't much matter.  And similarly, we now consider the number $x_0$ to be fixed.

\section{The key theorem and its proof }
\label{sec:key_theorem}

To summarize: we have noted that the functions $\lambda$ and $\mu$ defined above
must be bounded and can be taken to be non-increasing and in fact must both have
the limit 0 at $\infty$.  Based on this, we fixed numbers $m_0$, $M_0$, and
$x_0$ as above so that (\ref{eq1}) and (\ref{eq2}) are satisfied.  We are now
ready to state and prove our key theorem:

\begin{Theorem}
  Under the conditions stated in the paragraph just above, the Akra-Bazzi bounds
  (\ref{akra_bazzi_result}) for solutions $T$ of the recursion
  (\ref{general_recursive_formula}) hold for all $h\in H(\lambda, \mu)$ iff both
  \begin{equation}
    \label{int_beta_over_t_converges}
    \int_1^\infty\frac{\lambda(t)}{t}\dt < \infty \text{ \qquad and \qquad }
    \int_1^\infty\frac{\mu(t)}{t}\dt < \infty
  \end{equation}
\end{Theorem}
\begin{Remark}
  All that really matters, since $\lambda$ and $\mu$ are bounded, is that the
  integral converges at $\infty$.  So when we write expressions for such
  functions, the understanding is that these expressions are valid for large
  values of $t$, and we are not concerned with possible unboundedness of these
  expressions otherwise.

  For example, we might take for any fixed $\alpha > 1$,
  \begin{equation*}
    \lambda(x) = \mu(x) = \frac{1}{\log^\alpha x}
  \end{equation*}
  These are the functions that Leighton shows suffice.  And Leighton also
  observes that we cannot take $\alpha = 1$; this will follow as well from the
  converse result below at the end of the proof.
\end{Remark}

\begin{Proof}\ 
  In what follows we will perform an induction using the multiplicative
  structure of the positive reals.  We will use the sequence $I_n$ of intervals
  that partition $[1,\infty)$ and which are defined as follows:
  \begin{alignat*}{2}
    I_0 &= [1,x_0m)\\
    I_j &= [x_0m_0^j, x_0m_0^{j+1}) &\quad& \text{for $j > 0$}\\
  \end{alignat*}
  \begin{enumerate}
  \item First we deal with sufficiency: we will show that
    (\ref{int_beta_over_t_converges}) suffices for the Akra-Bazzi bounds.

    Now the Akra-Bazzi bounds really consist of a lower bound and an upper bound,
    and we will deal with them separately.  And starting out, we will not even
    assume (\ref{int_beta_over_t_converges}), but simply look at what we can say
    about such bounds a priori.

    So first we deal with the lower bound.  We will produce a decreasing sequence
    of numbers $\gamma_j > 0$ such that
    \begin{equation}
      \label{eq_gamma}
      T(x) \geq \gamma_j x^p\Bigl(1 + \int_{x_0}^x t^{-p}g(t)\frac{dt}{t}\Bigr)
      \text{ \qquad for } x\in I_j
    \end{equation}

    Remembering that $1 < m < M < x_0$, then since $T$ is bounded away from 0 on
    $[1,x_0m]$, we know that $\gamma_0$ exists.  And by making $\gamma_0$ smaller
    if necessary, we can ensure that $\gamma_0D < 2^{-p}$.  We do that.
    $\gamma_0$ is henceforth fixed.

    We will now show by induction on $j$, using the recursion
    (\ref{general_recursive_formula}), that $\gamma_j$ exists for all $j \geq 0$,
    and that the sequence \set[[\gamma_j]] is very well behaved.

    Suppose we have found a decreasing sequence of numbers
    $\set[[\gamma_k:0\leq k < j]]$ satisfying (\ref{eq_gamma}).  We will show then
    how $\gamma_j$ can be defined satisfying the same constraint.  (And of course
    since the sequence is decreasing we will also have $\gamma_kD < 2^{-p}$ for
    each $k$.)

    So we assume that (\ref{eq_gamma}) holds for $j-1$.

    We know that (\ref{general_recursive_formula}) holds in any case:
    \begin{equation*}
      T(x) = \int_m^MT\Bigl(\frac{x}{t} + h(x,t)\Bigr)\dmu(t) + g(x)
    \end{equation*}
    and we are considering values of $x$ in the range
    $x_0m_0^j \leq x < x_0m_0^{j+1}$.  We know by (\ref{eq1}) that
    \begin{equation*}
      1 < \frac{x}{M_0} \leq \frac{x}{t} + h(x,t) \leq \frac{x}{m_0} < x_0m_0^j
    \end{equation*}
    so in any case, $\frac{x}{t} + h(x,t)$ is in one of the intervals we have
    already considered inductively, and so $T(x)$ is bounded below as desired by
    one of the $\gamma_j$ already defined.  And since all of those $\gamma_j$
    values are decreasing, they are all $\geq\gamma_{j-1}$.  Thus we have

    \begin{align*}
      T(x) & = \int_m^M T\biggl(\frac{x}{t} + h(x,t)\biggr)\dmu(t) + g(x) \\
           &\geq \gamma_{j-1}\int_m^M \biggl(\frac{x}{t} + h(x,t)\biggr)^p
             \biggl(1 + \int_{x_0}^{\frac{x}{t} + h(x,t)}s^{-p} g(s)\frac{ds}{s}\biggr)\dmu(t)
             + g(x) \\
           &= \gamma_{j-1}\int_m^M t^{-p}x^p\biggl(1 + t\frac{h(x,t)}{x}\biggr)^p
             \biggl(1 + \int_{x_0}^x s^{-p} g(s)\frac{ds}{s}
             - \int_{\frac{x}{t} + h(x,t)}^x s^{-p}g(s)\frac{ds}{s}\biggr)\dmu(t) + g(x)
    \end{align*}
    Now since $\lambda$ is decreasing, we have, using (\ref{eq2}), that
    \begin{equation*}
      1 - M\lambda(x_0m_0^j)
      \leq 1 - M\lambda(x)
      \leq 1 + t\frac{h(x,t)}{x}
    \end{equation*}
    so this becomes
    \begin{align*}
      T(x) \geq&\qquad \gamma_{j-1}\bigl(1 - M\lambda(x_0m_0^j)\bigr)^p x^p
                 \biggl(1 + \int_{x_0}^x s^{-p} g(s)\frac{ds}{s}\biggr) \int_m^M
                 t^{-p}\dmu(t) \\
               &\quad - \gamma_{j-1}\bigl(1 - M\lambda(x_0m_0^j)\bigr)^p x^p
                 \biggl(\int_{\frac{x}{t} + h(x,t)}^x
                 s^{-p}g(s)\frac{ds}{s}\biggr) \int_m^M t^{-p}\dmu(t) \\
               &\quad + g(x)
    \end{align*}

    and then, using (\ref{formula_for_p}), (\ref{bounds_on_g_integral}), and the
    fact that $\gamma_{j-1}D\cdot 2^p \leq \gamma_0D\cdot 2^p < 1$, this
    becomes
    \begin{align*}
      T(x)
      &\geq \gamma_{j-1}\bigl(1 -  M\lambda(x_0m_0^j)\bigr)^p x^p
        \biggl(1 + \int_{x_0}^x s^{-p} g(s)\frac{ds}{s}\biggr)
        v    - \gamma_{j-1}D\cdot 2^p g(x) + g(x)\\
      &\geq \gamma_{j-1}\bigl(1 -  M\lambda(x_0m_0^j)\bigr)^p x^p
        \biggl(1 + \int_{x_0}^x s^{-p} g(s)\frac{ds}{s}\biggr)
    \end{align*}
    and thus we may take
    \begin{equation*}
      \gamma_j = \bigl(1 - M\lambda(x_0m_0^j)\bigr)^p \gamma_{j-1}
    \end{equation*}
    As promised, the sequence \set[[\gamma_j]] that this produces is decreasing.

    Now we deal with the upper bound.  We will produce an increasing sequence of
    numbers $\Gamma_j > 0$ such that
    \begin{equation}
      \label{eq_Gamma}
      T(x) \leq \Gamma_j x^p\Bigl(1 + \int_{x_0}^x t^{-p}g(t)\frac{dt}{t}\Bigr)
      \text{ \qquad for } x\in I_j
    \end{equation}
    Since $T$ is bounded on $[1,m_0x]$, we know that $\Gamma_0$ exists, and by
    making it larger if necessary we can ensure that $\Gamma_0C > 2^p$.

    Suppose then that we have found an increasing sequence of numbers
    $\set[[\Gamma_k:0\leq k < j]]$ satisfying (\ref{eq_Gamma}).  We will show how
    $\Gamma_j$ can be defined satisfying the same constraint.  (And of course
    since the sequence is increasing we also have $\Gamma_kC > 2^p$ for each $k$.)
    
    We have $x_0m_0^j \leq x < x_0m_0^{j+1}$.  Using (\ref{eq1}) again we have

    \begin{align*}
      T(x) & = \int_m^M T\biggl(\frac{x}{t} + h(x,t)\biggr)\dmu(t) + g(x) \\
           &\leq \Gamma_{j-1}\int_m^M\biggl(\frac{x}{t} + h(x,t)\biggr)^p
             \biggl(1 + \int_{x_0}^{\frac{x}{t} + h(x,t)}s^{-p} g(s)\frac{ds}{s}\biggr)\dmu(t)
             + g(x) \\
           &= \Gamma_{j-1}\int_m^M t^{-p}x^p\biggl(1 + t\frac{h(x,t)}{x}\biggr)^p
             \biggl(1 + \int_{x_0}^x s^{-p} g(s)\frac{ds}{s}
             - \int_{\frac{x}{t} + h(x,t)}^x s^{-p}g(s)\frac{ds}{s}\biggr)\dmu(t) + g(x)
    \end{align*}

    Now since $\mu$ is decreasing, we have, using (\ref{eq2}), that
    \begin{equation*}
      1 + t\frac{h(x,t)}{x}
      \leq 1 + M\mu(x)
      \leq 1 + M\mu(x_0m_0^j)
    \end{equation*}
    so this becomes
    \begin{align*}
      T(x) \leq&\qquad \Gamma_{j-1e}\bigl(1 + M\mu(x_0m_0^j)\bigr)^px^p
                 \biggl(1 + \int_{x_0}^xs^{-p}g(s)\frac{\ds}{s}\biggr)
                 \int_m^Mt^{-p}\dmu(t)  \\
               &\quad -\Gamma_{j-1}\bigl(1 + M\mu(x_0m_0^j)\bigr)^px^p
                 \biggl(\int_{\frac{x}{t} +
                 h(x,t)}^xs^{-p}g(s)\frac{ds}{s}\biggr)\int_m^Mt^{-p}\dmu(t) \\
               &\quad + g(x)
    \end{align*}

    and so continuing, using (\ref{bounds_on_g_integral}) and
    (\ref{formula_for_p}), we have
    \begin{align*}
      T(x)
      &\leq \Gamma_{j-1}\bigl(1 + M\mu(x_0m_0^j)\bigr)^p x^p
        \biggl(1 + \int_{x_0}^x t^{-p} g(t)\frac{dt}{t}\biggr)
        -\Gamma_{j-1}C\cdot\Bigl(\frac{1}{2}\Bigr)^pg(x) + g(x) \\
      &\leq \Gamma_{j-1}\bigl(1 + M\mu(x_0m_0^j)\bigr)^p x^p
        \biggl(1 + \int_{x_0}^x t^{-p} g(t)\frac{dt}{t}\biggr)
    \end{align*}
    and so we may take
    \begin{equation*}
      \Gamma_j = \bigl(1 + M\mu(x_0 m_0^j)\bigr)^p\Gamma_{j-1}
    \end{equation*}

    Now for (\ref{akra_bazzi_result}) to hold, it suffices that both $\gamma_j$
    and $\Gamma_j$ have finite non-zero limits.  And this is true iff the products
    \begin{equation*}
      \prod\bigl(1 - M\lambda(x_0m_0^j)\bigr)^p
      \quad\text{ and }\quad
      \prod\bigl(1 + M\mu(x_0m_0^j)\bigr)^p
    \end{equation*}
    both converge.  This is equivalent to the convergence of the sums
    \begin{equation*}
      \sum\lambda(x_0 m_0^j) \text{ \qquad and \qquad } \sum\mu(x_0 m_0^j)
    \end{equation*}
    Since $\lambda$ and $\mu$ are decreasing, this in turn is equivalent to the
    convergence of the integrals
    \begin{equation*}
      \int_1^\infty \frac{\lambda(t)}{t}\dt
      \text{ \qquad and \qquad } \int_1^\infty \frac{\mu(t)}{t}\dt
    \end{equation*}

    Thus we have shown that the convergence of these integrals suffices for
    (\ref{akra_bazzi_result}).

  \item Now we prove the converse.  Suppose that the integral for $\mu$
    diverges.  We just need to provide a counterexample.

    We will take a very simple function $T$, defined in the following way: we
    take $g\equiv 0$, and pick numbers $a>1$, $b=3$, and $p$ with $a = b^p$.  We
    take $h(x,3) = x\mu(x)$.  We take $m_0 = 2$.

    We then set $T(x) = 1$ for
    $1\leq x \leq x_0$ and we define $T(x)$ inductively for $x > x_0$ by
    \begin{equation*}
      T(x) = aT\Bigl(\frac{x}{b} + h(x, b)\Bigr)
    \end{equation*}

    This works because for $\abs{h(x,b)} < \frac{x}{12}$
    \begin{equation*}
      \frac{x}{4} \leq \frac{x}{3} + h(x,3) < \frac{x}{2}
    \end{equation*}
    so we only have to make $x_0$ large enough so that this is true for all $x > x_0$.

    We will next produce inductively a sequence of positive numbers
    $\set[[\theta_j: 0 \leq j < \infty]]$ such that for
    \begin{equation*}
      x_0m_0^j \leq x < x_0m_0^{j+1}
    \end{equation*}
    we have the lower bound
    \begin{equation*}
      T(x) \geq \theta_jx^p
    \end{equation*}
    Certainly $\theta_0$ exists.  Once we have produced $\theta_{j-1}$, we can
    (remembering that $b = 3$) compute for $x_0m_0^j \leq x < x_0m_0^{j+1}$ that
    \begin{align*}
      T(x) &= aT\Bigl(\frac{x}{b} + h(x, b)\Bigr) \\
           &\geq \theta_{j-1}a\Bigl(\frac{x}{b} + h(x, b)\Bigr)^p \\
           &= \theta_{j-1}ab^{-p}x^p\bigl(1 + b\mu(x)\bigr)^p \\
           &\geq \theta_{j-1}\bigl(1 + b\mu(x_0m_0^{j+1}\bigr)^p x^p
    \end{align*}
    so by induction we have 
    \begin{equation*}
      \theta_j \geq \theta_0\prod_{k=1}^{j+1}\bigl(1 + b\mu(x_0m_0^k\bigr)^p
    \end{equation*}
    which increases to $\infty$, contradicting (\ref{akra_bazzi_result}).  This
    subsumes Leighton's example of $h(x) = x/\log x$.

    An entirely similar computation works if the integral for $\lambda$
    diverges.  The proof is complete.
  \end{enumerate}
\end{Proof}

\appendix

\section{Functions of  O-regular variation: some equivalent properties}
\label{sec:OR}

All the material in this appendix is taken---not necessarily in this
order---from Chapter 2 of \ocite{Bingham-Goldie-Teugels1987}.

Suppose $f:\R^+\to\R^+$ is is measurable.  We define
\begin{align*}
  f^*(\lambda) &\defeq \limsup_{x\to\infty}\frac{f(\lambda x)}{f(x)} \\
  f_*(\lambda) &\defeq \liminf_{x\to\infty}\frac{f(\lambda x)}{f(x)}
\end{align*}
Using weak restrictions on $f^*$ and $f_*$, we can deduce some remarkable bounds
on the growth of $f$.

First, it can be shown that if $f^*(\lambda) < \infty$ on a set in $\R^*$ of
positive measure, then in fact there is a number $a_0 \geq 1$ such that
$f^*(\lambda) < \infty$ for all $\lambda \geq a_0$.

Further, it can be shown that if in addition there is at least one value
$\lambda_0 \geq 2a_0$ such that $f_*(\lambda_0) > -\infty$, then in fact
\begin{equation*}
  -\infty < f_*(\lambda) \leq f^*(\lambda) < \infty
\end{equation*}
for all $\lambda > 0$.

Such functions are said to be in the class \bfOregR of
\textbf{O-regularly varying functions}.

There are a number of equivalent characterizations of this class.  Here are
some, starting with the two just stated:

\renewcommand\labelenumi{(\arabic{enumi})}
\renewcommand\theenumi\labelenumi
\begin{enumerate}
\item\label{mult:OR_f_star_1} $f^*(\lambda) < \infty$ on a set of measure $>
  0$ and $f_*(\lambda) > 0$ for at least one $\lambda > 2a_0$;
\item\label{mult:OR_f_star_2} for all $\lambda \geq 1$,
  \begin{equation*}
    0 < f_*(\lambda) \leq f^*(\lambda) < \infty
  \end{equation*}
\item\label{mult:OR_multiplicative_polynomial_growth} there are numbers
  \begin{equation*}
    -\infty < \alpha < \beta < \infty
  \end{equation*}
  and positive constants $A$, $B$, and $X$ such that for $X \leq x \leq y$,
  \begin{equation*}
    B\Bigl(\frac{y}{x}\Bigr)^\beta \leq \frac{f(y)}{f(x)} \leq
    A\Bigl(\frac{y}{x}\Bigr)^\alpha
  \end{equation*}
\item\label{mult:OR_local_polynomial_growth} for each $c_0 \geq 1$ and
  $c_1 \geq 1$ with at least one of $c_0$ and $c_1$ strictly greater than 1,
  there are positive constants $C_0 \leq C_1$ and $X$ such that for all $x > X$ and
  \begin{equation*}
    \frac{x}{c_0} \leq y \leq c_1x
  \end{equation*}
  we have
  \begin{equation*}
    C_0 \leq \frac{f(y)}{f(x)} \leq C _1
  \end{equation*}
\item\label{mult:OR_local_polynomial_growth_restricted} for some number $c_0 \geq 1$ and
  some number $c_1 \geq 1$ with at least one of $c_0$ and $c_1$ strictly greater than 1,
  there are positive constants $C_0 \leq C_1$ and $X$ such that for all $x > X$ and
  \begin{equation*}
    \frac{x}{c_0} \leq y \leq c_1x
  \end{equation*}
  we have
  \begin{equation*}
    C_0 \leq \frac{f(y)}{f(x)} \leq C _1
  \end{equation*}
\end{enumerate}

The implications
\begin{equation*}
  \ref{mult:OR_multiplicative_polynomial_growth}
  \implies\ref{mult:OR_local_polynomial_growth}
  \implies\ref{mult:OR_local_polynomial_growth_restricted}
  \implies\ref{mult:OR_f_star_1}
\end{equation*}
are straightforward\footnote{And
  $\ref{mult:OR_local_polynomial_growth_restricted}\implies\ref{mult:OR_local_polynomial_growth}$
  is also elementary.  Lemma~\ref{lemma:stretching_interval} gives us all we need
  here for this}.  The remaining implications
\begin{equation*}
  \ref{mult:OR_f_star_1} \implies \ref{mult:OR_f_star_2}
  \implies \ref{mult:OR_multiplicative_polynomial_growth}
\end{equation*}
are much more difficult and are proved in \ocite{Bingham-Goldie-Teugels1987},
Chapter 2.

Just as a remark, \ref{mult:OR_local_polynomial_growth} and
\ref{mult:OR_local_polynomial_growth_restricted} are of course true without
the constraint that at least one of $c_0$ and $c_1$ is strictly greater than 1.
But this constraint is necessary in order to show that either of these
conditions implies \ref{mult:OR_f_star_1}.

Statements \ref{mult:OR_local_polynomial_growth} and
\ref{mult:OR_local_polynomial_growth_restricted} are those used in this report.

\bibliography{divide_and_conquer.bib}

\begin{bibdiv}
\begin{biblist}

\bib{AkraBazzi1998}{article}{
      author={Akra, Mohamed},
      author={Bazzi, Louay},
       title={{On the Solution of Linear Recurrence Equations}},
        date={1998},
     journal={Computational Optimization and Applications},
      volume={10},
      number={2},
       pages={195\ndash 210},
}

\bib{BentleyHakenSaxe1980}{article}{
      author={Bentley, Jon~Louis},
      author={Haken, Dorothea},
      author={Saxe, James~B.},
       title={{A General Method for Solving Divide-and-Conquer Recurrences}},
        date={1980},
     journal={SIGACT News},
      volume={12},
      number={3},
       pages={36\ndash 44},
}

\bib{Bingham-Goldie-Teugels1987}{book}{
      author={Bingham, N.~H.},
      author={Goldie, C.~M.},
      author={Teugels, J.~L.},
       title={{Regular Variation}},
   publisher={Cambridge Univerisity Press},
     address={Cambridge (U. K.)},
        date={1987},
        note={[QA331.5.B54]},
}

\bib{CormenLeisersonRivestStein2009}{book}{
      author={Cormen, Thomas~H.},
      author={Leiserson, Charles~E.},
      author={Rivest, Ronald~L.},
      author={Stein, Clifford},
       title={{Introduction to Algorithms}},
     edition={Third Edition},
   publisher={MIT Press},
     address={Cambridge, Massachusetts and London, England},
        date={2009},
        note={[QA76.6.I5858 2009]},
}

\bib{Kao1997}{article}{
      author={Kao, Ming-Yang},
       title={{Multiple-Size Divide-and-Conquer Recurrences}},
        date={1997},
     journal={SIGACT News},
      volume={28},
      number={2},
       pages={67\ndash 69},
        note={Earlier version in \emph{Proceedings of the International
  Conference on Algorithms, the 1996 International Computer Symposium}, pages
  159--161, National Sun Yat-Sen University, Kaohsiung, Taiwan, Republic of
  China, 1996},
}

\bib{Leighton1996}{misc}{
      author={Leighton, Tom},
       title={{Notes on Better Master Theorems for Divide-and-Conquer
  Recurrences}},
        date={1996},
        note={Lecture Notes; available at
  \url{https://courses.csail.mit.edu/6.046/spring04/handouts/akrabazzi.pdf}},
}

\bib{Wang_Fu1996}{article}{
      author={Wang, Xiaodong},
      author={Fu, Qingxiang},
       title={A frame for general divide-and-conquer recurrences},
        date={1996},
     journal={Information Processing Letters},
      volume={59},
       pages={45\ndash 61},
}

\end{biblist}
\end{bibdiv}
\end{document}